\documentclass[12pt]{article}

\usepackage{amsxtra,amssymb,amsthm,amsmath,latexsym}
\usepackage{graphicx}
\usepackage{epsfig}
\usepackage{epstopdf}
\usepackage{booktabs}
\usepackage{multirow}
\usepackage{afterpage}
\usepackage{setspace}
\usepackage{microtype}
\usepackage{calc,fourier}
\usepackage{cite}
\usepackage[T1]{fontenc}
\usepackage[hidelinks]{hyperref}
\usepackage{float}
\usepackage[small]{caption}

\numberwithin{equation}{section}

\newcommand{\bee}{\begin{equation*}}
\newcommand{\eee}{\end{equation*}}
\newcommand{\be}{\begin{equation}}
\newcommand{\ee}{\end{equation}}
\newcommand{\ba}{\begin{align}}
\newcommand{\ea}{\end{align}}

\newcommand{\RRR}{\mathbb{R}^3}

\title{Numerical Method for Solving Electromagnetic Wave Scattering by One and Many Small Perfectly Conducting Bodies}
\author{N. T. Tran\footnote{Mailing address:  Mathematics Department, 138 Cardwell Hall, Manhattan, KS 66506} \\
\small Department of Mathematics\\[-0.8ex]
\small Kansas State University, Manhattan, KS 66506-2602, USA\\
\small \texttt{nhantran@ksu.edu}}

\date{}
\begin{document}
\maketitle

\begin{abstract}
In this paper, we investigate the problem of electromagnetic (EM) wave scattering by one and many small perfectly conducting bodies and present a numerical method for solving it. For the case of one body, the problem is solved for a body of arbitrary shape, using the corresponding boundary integral equation. For the case of many bodies, the problem is solved asymptotically under the physical assumptions $a\ll d \ll \lambda$, where $a$ is the characteristic size of the bodies, $d$ is the minimal distance between neighboring bodies, $\lambda=2\pi/k$ is the wave length and $k$ is the wave number. Numerical results for the cases of one and many small bodies are presented. Error analysis for the numerical method are also provided.
\end{abstract}



\section{Introduction} \label{sec0}
Many real-world electromagnetic (EM) problems like EM wave scattering, EM radiation, etc \cite{Tsang2004}, cannot be solved analytically and exactly to get a solution in a closed form. Thus, numerical methods have been developed to tackle these problems approximately. Computational Electromagnetics (CEM) has evolved enormously in the past decades to a point that its methods can solve EM problems with extreme accuracy. These methods can be classified into two categories: Integral Equation (IE) method and Differential Equation (DE) method. Typical IE methods include: Method of Moment (MoM) developed by Roger F. Harrington (1968) \cite{harrington1996field}, Fast Multipole Method (FMM) first introduced by Greengard and Rokhlin (1987) \cite{GR1987} and then applied to EM by Engheta et al (1992) \cite{engheta1985fast}, Partial Element Equivalent Circuit (PEEC) method \cite{ruehli1974equivalent}, and Discrete Dipole Approximation \cite{devoe1964optical}. Typical DE methods are: Finite Difference Time Domain (FDTD) developed by Kane Yee (1966) \cite{yee1966numerical}, Finite Element Method (FEM) \cite{zienkiewicz1977finite}, Finite Integration Technique (FIT) proposed by Thomas Weiland (1977) \cite{weiland2001discrete}, Pseudospectral Time Domain (PSTD) \cite{liu1997pstd}, Pseudospectral Spatial Domain (PSSD) \cite{tyrrell2005pseudospectral}, and Transmission Line Matrix (TLM) \cite{hoefer1985transmission}. Among these methods, FDTD has emerged as one of the most popular techniques for solving EM problems due to its simplicity and ability to provide animated display of the EM field. However, FDTD requires the entire computational domain be gridded \cite{stavroulakis2013biological}, that results in very long solution times. Furthermore, as a DE method, it does not take into account the radiation condition in exact sense \cite{kunz1993finite, buchanan1993simulation}, which leads to certain error in the solution. On the other hand, spurious solutions might exist in DE methods \cite{mur2002causes, lu2004elimination, zhao2007spurious}. Most importantly, most of DE methods are not suitable if the number of bodies is very large.

In \cite{R620,R635,R652,andriychuk2012electromagnetic,ramm2014application,ramm2014calculation}, A. G. Ramm has developed a theory of EM wave scattering by many small perfectly conducting and impedance bodies. In this theory, the EM wave scattering problem is solved asymptotically under the physical assumptions: $a\ll d \ll \lambda$, where $a$ is the characteristic size of the bodies, $d$ is the minimal distance between neighboring bodies, $\lambda=2\pi/k$ is the wave length and $k$ is the wave number. In \cite{TEMImpedance}, a numerical method is developed for solving  EM wave scattering by many small impedance bodies. In this paper, the problem of EM wave scattering by one and many small perfectly conducting bodies is considered. A numerical method for solving this problem asymptotically based on the above theory is presented. For the case of one body, the problem is solved for a body of arbitrary shape, using the corresponding boundary integral equation. For the case of many small bodies, the problem is solved under the basic assumptions $a\ll d \ll \lambda$ and the assumption about the distribution of the small bodies
\be
    \mathcal{N}(\Delta)=\frac{1}{a^3} \int_\Delta N(x)dx[1+o(1)], \quad a\to 0,
\ee
in which $\Delta$ is an arbitrary open subset of the domain $\Omega$ that contains all the small bodies, $\mathcal{N}(\Delta)$ is the number of the small bodies in $\Delta$,  and $N(x)$ is the distribution function of the bodies
\be
    N(x)\ge 0, \quad N(x) \in C(\Omega).
\ee

In Sections \ref{sec1} and \ref{sec2}, the theory of EM wave scattering by one and many small perfectly conducting bodies is presented. The numerical methods for solving these problems are also  described in details. Furthermore, error analysis for the numerical methods of solving the EM scattering problem are also provided. In Section \ref{sec3} these methods are tested and numerical results are discussed. 

\section{EM wave scattering by one perfectly conducting body} \label{sec1}
Let $D$ be a bounded perfectly conducting body, $a=\frac{1}{2}$diam$D$, $S$ be its $C^2$-smooth boundary,  and $D':=\RRR \setminus D$. Let $\epsilon$ and $\mu$ be the dielectric permittivity and magnetic permeability constants of the medium in $D'$. Let $E$ and $H$ denote the electric and magnetic fields, respectively, $E_0$ be the incident field and $v_E$ be the scattered field. The problem of electromagnetic wave scattering by one perfectly conducting body can be stated as follows
\begin{align}
            &\nabla \times E=i \omega \mu H, \quad \text{in }D':=\RRR \setminus D, \label{eq1.1.1} \\
            &\nabla \times H=-i \omega \epsilon E, \quad \text{in }D', \label{eq1.1.2} \\
            &[N,[E, N]]=0, \quad \text{on } S:=\partial D, \label{eq1.1.3} \\
            &E=E_0+v_E,  \label{eq1.1.4} \\
            &E_0=\mathcal{E} e^{ik\alpha\cdot x}, \quad \mathcal{E} \cdot \alpha=0,
            \quad \alpha \in S^2, \label{eq1.1.5} \\
            &\frac{\partial v_E}{\partial r}-ikv_E=o\left(\frac{1}{r}\right), \quad r:=|x|\to \infty, \label{eq1.1.6}
\end{align}
where $\omega>0$ is the frequency, $k=2\pi/\lambda=\omega \sqrt{\epsilon \mu}$ is the wave number, $ka \ll 1$, $\lambda$ is the wave length, $\mathcal{E}$ is a constant vector, and $\alpha$ is a unit vector that indicates the direction of the incident wave $E_0$. This incident wave satisfies the relation $\nabla\cdot E_0=0$. The scattered field $v_E$ satisfies the radiation condition \eqref{eq1.1.6}. Here, $N$ is the unit normal vector to the surface $S$, pointing out of $D$. By $[\cdot,\cdot]$ the vector product is denoted and $\alpha \cdot x$ is the scalar product of two vectors.

The solution to problem \eqref{eq1.1.1}-\eqref{eq1.1.6} can be found in the form
\be \label{eq1.3.2}
    E(x)=E_0(x)+\nabla \times \int_S g(x,t) J(t)dt, \quad g(x,t):=\frac{e^{ik|x-t|}}{4\pi |x-t|},
\ee
see \cite{R620}.
Here, $E$ is a vector in $\RRR$ and $\nabla \times E$ is a pseudo-vector, that is a vector-like object which changes sign under reflection of its coordinate axes. $E_0$ is the incident plane wave defined in \eqref{eq1.1.5} and $J$ is an unknown pseudo-vector that is to be found. $J$ is assumed to be tangential to $S$ and continuous. $J$ can be found by applying the boundary condition \eqref{eq1.1.3}, or equivalently $[N,E]=0$, to \eqref{eq1.3.2} and solving the resulting boundary integral equation
\be \label{eq1.3.2b}
    \frac{J}{2}+AJ:=\frac{J(s)}{2}+\int_S [N_s,[\nabla_s g(s,t),J(t)]]dt=-[N_s,E_0],
\ee
or, equivalently
\be \label{eq1.3.2c}
    (I+2A)J=F,
\ee
where $F:=-2[N_s,E_0]$. Equation \eqref{eq1.3.2c} is of Fredholm type since $A$ is compact, see \cite{R652}.

Once we have $J$, $E$ can be computed by formula \eqref{eq1.3.2} and $H$ can be found by the formula
\be \label{eq1.3.3}
    H=\frac{\nabla\times E}{i\omega \mu}.
\ee

If $D$ is sufficiently small, then equation \eqref{eq1.3.2c} is uniquely solvable in $C(S)$ and its solution $J$ is tangential to $S$, see \cite{R620}. The asymptotic formula for $E$ when the radius $a$ of the body $D$ tends to zero can be derived as follows, see \cite{R620}. Rewrite equation \eqref{eq1.3.2} as
\be \label{eq1.3.4}
    E(x)=E_0(x)+[\nabla g(x,x_1), Q]+ \nabla\times \int_{S} [g(x,t)-g(x,x_1)] J(t)dt,
\ee
where $x_1 \in D$, an arbitrary point inside the small body $D$, and
\be \label{eq1.3.5}
    Q:=\int_S J(t)dt.
\ee
Since
\begin{align}
	&|\nabla g(x,x_1)]| = O\left(\frac{k}{d}+\frac{1}{d^2}\right), \quad d=|x-x_1|,  \\
	&|g(x,t)-g(x,x_1)| = O\left(\left(\frac{k}{d}+\frac{1}{d^2}\right)a\right), \quad a=|t-x_1|, \quad \text{ and} \\
    &|\nabla[g(x,t)-g(x,x_1)]| = O\left(\frac{ak^2}{d}+\frac{ak}{d^2}+\frac{a}{d^3}\right),
\end{align}
the second term in \eqref{eq1.3.4} is much greater than the last term
\be \label{eq1.3.6}
    \left|[\nabla g(x,x_1), Q]\right| \gg \left|\nabla\times \int_{S} [g(x,t)-g(x,x_1)] J(t)dt\right|, \quad a \to 0.
\ee
Then, the asymptotic formula for $E$ when $a$ tends to zero is
\be \label{eq1.3.8}
    E(x)=E_0(x)+[\nabla_x g(x,x_1),Q], \quad a \to 0,
\ee
where $|x-x_1|\gg a$, $x_1 \in D$.
Thus, when $D$ is sufficiently small, instead of finding $J$, we can just find one pseudovector $Q$.

The analytical formula for $Q$ is derived as follows, see \cite{R652}. By integrating both sides of \eqref{eq1.3.2b} over $S$, one gets
\be \label{eq1.3.9}
    \int_S \frac{J(s)}{2}ds +\int_S ds \int_S dt [N_s,[\nabla_s g(s,t),J(t)]]=-\int_S [N_s,E_0]ds.
\ee
This is equivalent to
\be \label{eq1.3.10}
    \frac{Q}{2}+\int_S dt \int_S ds \nabla_s g(s,t)N_s\cdot J(t)-\int_S dt J(t) \int_S ds \frac{\partial g(s,t)}{\partial N_s}=-\int_D \nabla\times E_0 dx.
\ee
When $a\to 0$, this equation becomes
\be \label{eq1.3.11}
    \frac{Q}{2}+e_p\int_S dt \int_S ds \frac{\partial g(s,t)}{\partial s_p}N_q(s) J_q(t)+\frac{1}{2}\int_S dt J(t) =-|D| \nabla\times E_0, \quad  1\le p,q\le 3,
\ee
where in the second term, summations over the repeated indices are understood, $e_p$, $1\le p\le 3$, are the orthogonal unit vectors in $\RRR$, $|D|$ is the volume of $D$, $|D|=c_D a^3$, and in the third term we use this estimate
\be \label{eq1.3.12}
    \int_S ds \frac{\partial g(s,t)}{\partial  N_s}\simeq \int_S ds \frac{\partial g_0(s,t)}{\partial  N_s}=-\frac{1}{2}, \quad g_0(s,t):=\frac{1}{4\pi|s-t|}.
\ee
Let
\be \label{eq1.3.13}
    \Gamma_{pq}(t):=\int_S ds \frac{\partial g(s,t)}{\partial  s_p}N_q(s),
\ee
then equation \eqref{eq1.3.11} can be rewritten as follows
\be \label{eq1.3.14}
    \frac{Q}{2}+e_p\int_S dt \Gamma_{pq}(t) J_q(t)+\frac{Q}{2}=-|D| \nabla\times E_0,
\ee
or
\be \label{eq1.3.15}
    Q+\Gamma Q=-|D| \nabla\times E_0,
\ee
where $\Gamma$ is a $3\times 3$ constant matrix and it is defined by
\be \label{eq1.3.16a}
    \Gamma Q=e_p\int_S dt \Gamma_{pq}(t) J_q(t),
\ee
in which summations are understood over the repeated indices. Thus, $Q$ can be written as
\be \label{eq1.3.16}
        Q=-|D|(I+\Gamma)^{-1} \nabla\times E_0, \quad a \to 0,
\ee
where $I:=I_3$, the $3\times 3$ identity matrix. This formula is asymptotically exact as $a \to 0$.

\subsection{Numerical method for solving EM wave scattering by one perfectly conducting spherical body}\label{sec3.1}
In this section, we consider the EM wave scattering problem by a small perfectly conducting spherical body. Instead of solving the problem \eqref{eq1.1.1}-\eqref{eq1.1.6} directly, we will solve its corresponding boundary integral equation \eqref{eq1.3.2b} for the unknown vector $J$
\be \label{eq3.0.1}
    \frac{J(s)}{2}+\int_S [N_s,[\nabla_s g(s,t),J(t)]]dt=-[N_s,E_0].
\ee
Then the solution $E$ to the EM wave scattering problem by one perfectly conducting body can be computed by
either the exact formula \eqref{eq1.3.2} or the asymptotic formula \eqref{eq1.3.8}.

Scattering by a sphere has been discussed in many papers, for example \cite{Mie1908} in which Mie solves the EM wave scattering problem by separation of variables. The EM field, scattered by a small body, is proportional to $O(a^3)$.

Suppose $S$ is a smooth surface of a spherical body. Let $S$ be partitioned into $P$ non-intersecting subdomains $S_{ij}, 1 \le i \le m_{\theta}, 1 \le j \le m_{\phi}$, using spherical coordinates, where  $m_{\theta}$ is the number of intervals of $\theta$ between $0$ and $2\pi$ and  $m_{\phi}$ defines the number of intervals of $\phi$ between $0$ and $\pi$. Then $P=m_{\theta}m_{\phi}+2$, which includes the two poles of the sphere. $m_{\theta}$ is defined in this way: $m_{\theta}=m_{\phi}+|\phi-\frac{\pi}{2}|6m_{\phi}$. This means the closer it is to the poles of the sphere, the more intervals for $\theta$ are used. Then the point $(\theta_i,\phi_j)$ in $S_{ij}$ is chosen as follows
\begin{align}
	&\theta_i = i\frac{2\pi}{m_{\theta}}, \quad  1 \le i \le m_{\theta}, \label{eq3.0.2} \\
	&\phi_j = j\frac{\pi}{m_{\phi}+1}, \quad 1 \le j \le m_{\phi}. \label{eq3.0.3}
\end{align}
Note that there are many different ways to distribute collocation points. However, the one that we describe here will guarantee convergence to the solution to \eqref{eq3.0.1} with fewer collocation points used from our experiment. Furthermore, one should be careful when choosing the distribution of collocation points on a sphere. If one chooses $\phi_j = j\frac{\pi}{m_{\phi}}, 1 \le j \le m_{\phi}$, then when $j= m_{\phi}$, $\phi_j=\pi$ and thus there is only one point for this $\phi$ regardless of the value of $\theta$ as shown in \eqref{eq3.0.4}. The position of a point in each $S_{ij}$ can be computed by
\be \label{eq3.0.4}
    (x,y,z)_{ij} = a(\cos\theta_i \sin\phi_j,\sin\theta_i \sin\phi_j,\cos\phi_j),
\ee
and the outward-pointing unit normal vector $N$ to $S$ at this point is
\be \label{eq3.0.5}
    N_{ij}=N(\theta_i,\phi_j)=(\cos\theta_i \sin\phi_j,\sin\theta_i \sin\phi_j,\cos\phi_j).
\ee
For a star-shaped body with a different shape, only the normal vector $N$ needs to be recomputed. Rewrite the integral equation \eqref{eq3.0.1} as
\be \label{eq3.0.6a}
    \frac{J(s)}{2}+\int_S \nabla_s g(s,t)N_s\cdot J(t)dt-\int_S \frac{\partial g(s,t)}{\partial N_s} J(t)dt=-[N_s,E_0].
\ee
This integral equation can be discretized as follows
\be \label{eq3.0.6}
    J(i)+2\sum_{j \neq i}^P [\nabla_s g(i,j)N_s(i)\cdot J(j)-J(j)\nabla_s g(i,j)\cdot N_s(i)]\Delta_{j} = F(i), \quad 1 \le i \le P,
\ee
in which by $i$ the point $(x_i,y_i,z_i)$ is denoted, $F(i):=-2[N_s,E_0](i)$, and $\Delta_j$ is the surface area of the subdomain $j$. This is a linear system with unknowns $J(i):=(X_i,Y_i,Z_i), 1 \le i \le P$. This linear system can be rewritten as follows
\begin{align}
    &X_i+\sum_{j\neq i}^P a_{ij}X_j+ b_{ij}Y_j+ c_{ij}Z_j = F_x(i), \label{eq3.0.7}\\
    &Y_i+\sum_{j\neq i}^P a'_{ij}X_j+ b'_{ij}Y_j+ c'_{ij}Z_j = F_y(i), \label{eq3.0.8}\\
    &Z_i+\sum_{j\neq i}^P a''_{ij}X_j+ b''_{ij}Y_j+ c''_{ij}Z_j = F_z(i), \label{eq3.0.9}
\end{align}
where by the subscripts $x,y,z$ the corresponding coordinates are denoted, e.g. $F(i)=(F_x,F_y,F_z)(i)$, and
\begin{align}
    &a_{ij}:=2[\nabla g(i,j)_x N_x(i)-\nabla g(i,j)\cdot N(i)]\Delta_j,\\
    &b_{ij}:=2\nabla g(i,j)_x N_y(i)\Delta_j,\\
    &c_{ij}:=2\nabla g(i,j)_x N_z(i)\Delta_j,
\end{align}
for $i\neq j$; when $i=j$: $a_{ii}=1, b_{ii}=0$, and $c_{ii}=0$,
\begin{align}
    &a'_{ij}:=2\nabla g(i,j)_y N_x(i)\Delta_j,\\
    &b'_{ij}:=2[\nabla g(i,j)_y N_y(i)-\nabla g(i,j)\cdot N(i)]\Delta_j,\\
    &c'_{ij}:=2\nabla g(i,j)_y N_z(i)\Delta_j,
\end{align}
for $i\neq j$; when $i=j$: $a'_{ii}=0, b'_{ii}=1$, and $c'_{ii}=0$,
\begin{align}
    &a''_{ij}:=2\nabla g(i,j)_z N_x(i)\Delta_j,\\
    &b''_{ij}:=2\nabla g(i,j)_z N_y(i)\Delta_j,\\
    &c''_{ij}:=2[\nabla g(i,j)_z N_z(i)-\nabla g(i,j)\cdot N(i)]\Delta_j,
\end{align}
for $i\neq j$; when $i=j$: $a''_{ii}=0, b''_{ii}=0$, and $c''_{ii}=1$.

\subsection{Error analysis}\label{sec3.1.1}
Recall the boundary integral equation \eqref{eq1.3.2b}
\be \label{eq3.0.10}
    \frac{J(s)}{2}+\int_S [N_s,[\nabla_s g(s,t),J(t)]]dt=-[N_s,E_0].
\ee
Integrate both sides of this equation over $S$ and get
\be \label{eq3.0.11}
    Q+\Gamma Q=-|D| \nabla\times E_0,
\ee
see Section \ref{sec1}. Once $J$ is found from solving \eqref{eq3.0.10}, $Q$ can be computed by $Q=\int_S J(t)dt$.
Then one can validate the values of $J$ and $Q$ by checking the following things
\begin{itemize}
\item Is $J$ tangential to $S$ as shown in Section \ref{sec1}? One needs to check $J(s)\cdot N_s$.
\item Is $Q=\int_S J(t)dt$ correct? The relative error of $Q$ can be computed as follows
\be \label{eq3.0.12}
    \text{Error} = \frac{|Q+\Gamma Q-RHS|}{|RHS|},
\ee
where $RHS:=-|D| \nabla\times E_0$. This will give the error of the numerical method for the case of one body.
\end{itemize}
Furthermore, one can also compare the value of the asymptotic $Q_a$ in formula \eqref{eq1.3.16} with the exact $Q_e$ defined in \eqref{eq1.3.5} by
\be \label{eq3.0.13}
    \text{Error} = \frac{|Q_e-Q_a|}{|Q_e|},
\ee
and check the difference between the asymptotic $E_a$ in \eqref{eq1.3.8} and the exact $E_e$ defined in \eqref{eq1.3.2} by computing this relative error
\be \label{eq3.0.14}
    \text{Error} = \frac{|E_e-E_a|}{|E_e|}.
\ee

\subsection{General method for solving EM wave scattering by one perfectly conducting body}\label{sec3.4}
In this section, we present a general method for solving the EM wave scattering problem by one perfectly conducting body, whose surface is parametrized by $f(u,v) = (x(u,v), y(u,v),$ $ z(u,v))$.

\begin{itemize}
\item Step 1: One needs to partition the surface of the body into $P$ non-intersecting subdomains. In each subdomain, choose a collocation point. The position of the collocation points can be computed using $f(u,v)=(x(u,v),y(u,v),z(u,v))$, see for example \eqref{eq3.0.2}-\eqref{eq3.0.4}.
\item Step 2: Find the unit normal vector $N$ of the surface from the function $f$.
\item Step 3: Solve the linear system \eqref{eq3.0.7}-\eqref{eq3.0.9} for $X_i,Y_i$, and $Z_i, 1 \le i \le P$. Then vector $J$ in the boundary integral equation \eqref{eq1.3.2b} is computed by $J(i):=(X_i,Y_i,Z_i)$ at the point $i$ on the surface.
\item Step 4: Compute the electric field $E$ using \eqref{eq1.3.2}.
\end{itemize}

\section{EM wave scattering by many small perfectly conducting bodies} \label{sec2}
Consider a bounded domain $\Omega$ containing $M$ small bodies $D_m$, $1\le m \le M$, and $S_m$ are their corresponding smooth boundaries. Let $D:=\bigcup_{m=1}^M D_m \subset \Omega$ and $D'$ be the complement of $D$ in $\RRR$. We assume that $S=\bigcup_{m=1}^M S_m$ is C$^2$-smooth. $\epsilon$ is  the dielectric permittivity constant and $\mu$ is the magnetic permeability constant of the medium. Let $E$ and $H$ denote the electric and magnetic fields, respectively. $E_0$ is the incident field and $v$ is the scattered field. The problem of electromagnetic wave scattering by many small perfectly conducting bodies  involves solving the following system
\begin{align}
            &\nabla \times E=i \omega \mu H, \quad \text{in }D':=\RRR \setminus D, \quad D:=\bigcup_{m=1}^M D_m, \label{eq2.1.1} \\
            &\nabla \times H=-i \omega \epsilon E, \quad \text{in }D', \label{eq2.1.2} \\
            &[N,[E, N]]=0, \quad \text{on } S,  \label{eq2.1.3} \\
            &E=E_0+v,  \label{eq2.1.4} \\
            &E_0=\mathcal{E} e^{ik\alpha\cdot x}, \quad \mathcal{E} \cdot \alpha=0,
            \quad \alpha \in S^2. \label{eq2.1.5}
\end{align}
where $v$ satisfies the radiation condition \eqref{eq1.1.6}, $\omega>0$ is the frequency, $k=2\pi/\lambda$ is the wave number, $ka \ll 1$, $a:=\frac{1}{2}\max_m \text{diam}D_m$, and $\alpha$ is a unit vector that indicates the direction of the incident wave $E_0$. Furthermore,
\be \label{eq2.1.7}
    \epsilon=\epsilon_0, \quad \mu=\mu_0 \quad \text{ in } \Omega':=\RRR \setminus \Omega.
\ee
Assume that the distribution of small bodies $D_m$, $1 \le m \le M$, in $\Omega$ satisfies the following formula
\be \label{eq2.1.8}
    \mathcal{N}(\Delta)=\frac{1}{a^3} \int_\Delta N(x)dx[1+o(1)], \quad a\to 0,
\ee
where $\mathcal{N}(\Delta)$ is the number of small bodies in $\Delta$, $\Delta$ is an arbitrary open subset of $\Omega$, and $N(x)$ is the distribution function
\be \label{eq2.1.9}
    N(x)\ge 0, \quad N(x) \in C(\Omega).
\ee

Note that $E$ solves this equation
\be \label{eq2.1.11}
    \nabla\times\nabla\times E=k^2 E, \quad k^2=\omega^2\epsilon \mu,
\ee
if $\mu=$const.
Once we have $E$, then $H$ can be found from this relation
\be \label{eq2.1.10}
    H=\frac{\nabla\times E}{i\omega\mu}.
\ee
From \eqref{eq2.1.10} and \eqref{eq2.1.11}, one can get \eqref{eq2.1.2}. Thus, we need to find only $E$ which satisfies the boundary condition \eqref{eq2.1.3}. It was proved in \cite{R620} that under the radiation condition and the assumptions $a \ll d \ll \lambda$, the problem \eqref{eq2.1.1}-\eqref{eq2.1.5} has a unique solution and its solution is of the form
\be \label{eq2.1.13}
    E(x)=E_0(x)+\sum_{m=1}^M \nabla\times \int_{S_m} g(x,t) J_m(t)dt,
\ee
 where $J_m$ are unknown continuous functions that can be found from the boundary condition.
 Let
\be \label{eq2.1.15}
    Q_m:=\int_{S_m}J_m(t)dt.
\ee
When $a \to 0$, the asymptotic solution for the electric field is given by
\be \label{eq2.1.22}
    E(x)=E_0(x)+\sum_{m=1}^M [\nabla g(x,x_m), Q_m], \quad a\to 0.
\ee
Therefore, instead of finding $J_m(t), \forall t \in S, 1 \le m \le M$, to get the solution $E$, one can just find $Q_m$. This allows one to solve the EM scattering problem with a very large number of small bodies which is impossible to do before. The analytic formula for $Q_m$ can be derived by using formula \eqref{eq1.3.16} and replacing $E_0$ in this formula by the effective field $E_{e}(x_m)$ acting on the m-th body
\be \label{eq2.1.24}
    Q_m=-|D_m|(I+\Gamma)^{-1} \nabla\times E_{e}(x_m), \quad 1 \le m \le M, \quad x_m \in D_m,
\ee
where the effective field acting on the m-th body is defined as
\be \label{eq2.1.23}
    E_e(x_m)=E_0(x_m)+\sum_{j\ne m}^M [\nabla g(x_m,x_j), Q_j]\quad 1 \le m \le M.
\ee
When $a \to 0$, the effective field $E_e(x)$ is asymptotically equal to the field $E(x)$ in \eqref{eq2.1.22} as proved in \cite{R620} and \cite{R652}.

Let $E_{em}:=E_{e}(x_m)$, where $x_m$ is a point in $D_m$. From \eqref{eq2.1.24}, and \eqref{eq2.1.23}, one gets
\be \label{eq2.1.28}
    E_{em}=E_{0m}-\sum_{j\ne m}^M [\nabla g(x_m,x_j), (I+\Gamma)^{-1}\nabla\times E_{ej}]|D_j|,\quad
       1 \le m \le M.
\ee

\subsection{Numerical method for solving EM wave scattering by many small perfectly conducting bodies}\label{sec3.5}
For finding the solution to EM wave scattering in the case of many small perfectly conducting bodies, we need to find $E_{em}$ in \eqref{eq2.1.28}. Apply the operator $(I+\Gamma)^{-1}\nabla \times$ to both sides of \eqref{eq2.1.28} and let $A_m:=(I+\Gamma)^{-1}\nabla \times E_{em}$. Then
\begin{align}
	A_m=A_{0m}-(I+\Gamma)^{-1}\sum_{j\ne m}^M |D_j|\left(\nabla_x\times [\nabla g(x,x_j),A_j] \right) |_{x=x_m}, \quad  1 \le m \le M, \label{eq3.5.1}
\end{align}
Solving this system yields $A_m$, for $1 \le m \le M$. Then $E$ can be computed by
\be \label{eq3.5.2}
    E(x)=E_0(x)+\sum_{m=1}^M [\nabla g(x,x_m), Q_m],
\ee
where
\be \label{eq3.5.2a}
    Q_m=-|D_m|A_m, \quad 1 \le m \le M.
\ee
Equation \eqref{eq3.5.1} can be rewritten as follows
\be \label{eq3.5.3}
    A_m=A_{0m}-\sum_{j\ne m}^M \tau[k^2 g(x_m,x_j)A_j+(A_j\cdot \nabla_x)\nabla g(x,x_j)|_{x=x_m}]|D_j|,
\ee
where $1 \le m \le M$, $\tau:=(I+\Gamma)^{-1}$, and $A_m$ are vectors in $\RRR$.

Let $A_i:=(X_i,Y_i,Z_i)$ then one can rewrite the system \eqref{eq3.5.3} as
\begin{align}
    &X_i+\sum_{j\neq i}^M a_{ij}X_j +  b_{ij}Y_j +  c_{ij}Z_j = F_x(i), \label{eq3.5.7}\\
    &Y_i+\sum_{j\neq i}^M a'_{ij}X_j +  b'_{ij}Y_j +  c'_{ij}Z_j = F_y(i), \label{eq3.5.8}\\
    &Z_i+\sum_{j\neq i}^M a''_{ij}X_j +  b''_{ij}Y_j +  c''_{ij}Z_j = F_z(i), \label{eq3.5.9}
\end{align}
in which by the subscripts $x,y,z$ the corresponding coordinates are denoted, e.g. $F(i)=(F_x,F_y,F_z)(i)$, where $F(i):=A_{0i}$ and
\begin{align}
    &a_{ij}:=[k^2 g(i,j) + \partial_{x} \nabla g(i,j)_x]|D_j|\tau(1,1),\\
    &b_{ij}:=\partial_{y} \nabla g(i,j)_x|D_j|\tau(1,1),\\
    &c_{ij}:=\partial_{z} \nabla g(i,j)_x|D_j|\tau(1,1),
\end{align}
for $i\neq j$, here $\tau(1,1)$ is the entry (1,1) of matrix $\tau$ in \eqref{eq3.5.3}; when $i=j$: $a_{ii}=1, b_{ii}=0$, and $c_{ii}=0$;
\begin{align}
    &a'_{ij}:=\partial_{x} \nabla g(i,j)_y|D_j|\tau(2,2),\\
    &b'_{ij}:=[k^2 g(i,j) + \partial_{y} \nabla g(i,j)_y]|D_j|\tau(2,2),\\
    &c'_{ij}:=\partial_{z} \nabla g(i,j)_y|D_j|\tau(2,2),
\end{align}
for $i\neq j$; when $i=j$: $a'_{ii}=0, b'_{ii}=1$, and $c'_{ii}=0$;
\begin{align}
    &a''_{ij}:=\partial_{x} \nabla g(i,j)_z|D_j|\tau(3,3),\\
    &b''_{ij}:=\partial_{y} \nabla g(i,j)_z|D_j|\tau(3,3),\\
    &c''_{ij}:=[k^2 g(i,j) + \partial_{z} \nabla g(i,j)_z]|D_j|\tau(3,3),
\end{align}
for $i\neq j$; when $i=j$: $a'_{ii}=0, b'_{ii}=0$, and $c'_{ii}=1$.

\subsection{Error analysis}
The error of the solution to the EM wave scattering problem by many small perfectly conducting bodies can be estimated as follows. From the solution $E$ of the electromagnetic scattering problem by many small bodies given in \eqref{eq2.1.13}
\be \label{eq3.5.48}
    E(x)=E_0(x)+\sum_{m=1}^M \nabla\times \int_{S_m} g(x,t) J_m(t)dt,
\ee
we can rewrite it as
\be \label{eq3.5.49}
    E(x)=E_0(x)+\sum_{m=1}^M [\nabla g(x,x_m), Q_m]+\sum_{m=1}^M \nabla\times \int_{S_m} [g(x,t)-g(x,x_m)] J_m(t)dt.
\ee
Comparing this with the asymptotic formula for $E$ when $a \to 0$ given in \eqref{eq2.1.22}
\be \label{eq3.5.50}
    E(x)=E_0(x)+\sum_{m=1}^M [\nabla g(x,x_m), Q_m],
\ee
we have the error of this asymptotic formula is
\be \label{eq3.5.51}
    \text{Error}=\left|\sum_{m=1}^M \nabla\times \int_{S_m} [g(x,t)-g(x,x_m)] J_m(t)dt\right| \sim \frac{1}{4\pi}\left(\frac{ak^2}{d}+\frac{ak}{d^2}+\frac{a}{d^3}\right)\sum_{m=1}^M|Q_m|,
\ee
where $d=\min_m|x-x_m|$ and
\be \label{eq3.5.52}
    Q_m=-|D_m|(I+\Gamma)^{-1} \nabla\times E_{e}(x_m), \quad 1 \le m \le M, \quad x_m \in D_m, \quad a \to 0,
\ee
because
\be \label{eq3.5.53}
    |\nabla[g(x,t)-g(x,x_m)]| = O\left(\frac{ak^2}{d}+\frac{ak}{d^2}+\frac{a}{d^3}\right), \quad a=\max_m|t-x_m|.
\ee

\section{Experiments} \label{sec3}

\subsection{EM wave scattering by one perfectly conducting spherical body}\label{sec3.1.2}
To illustrate the idea of the numerical method, we use the following physical parameters to solve the EM wave scattering problem by one small perfectly conducting sphere, i.e solving the linear system \eqref{eq3.0.7}-\eqref{eq3.0.9}
\begin{itemize}
     \item Speed of wave, $c=(3.0E+10)$ cm/sec.
     \item Frequency, $\omega=(5.0E+14)$ Hz.
     \item Wave number, $k = (1.05E+05)$ cm$^{-1}$.
     \item Wave length, $\lambda= (6.00E-05)$ cm.
     \item Direction of incident plane wave, $\alpha = (0,  1,  0)$.
     \item Magnetic permeability, $\mu = 1$.
     \item Vector $\mathcal{E} = (1,  0,  0)$.
     \item Incident field vector, $E_0$: $E_{0}(x)=\mathcal{E} e^{ik\alpha\cdot x}$.
     \item The body is a sphere of radius $a$, centered at the origin.
\end{itemize}
We use GMRES iterative method, see \cite{GMRES}, to solve the linear system \eqref{eq3.0.7}-\eqref{eq3.0.9}. For a spherical body, matrix $\Gamma$ in \eqref{eq1.3.15} can be computed analytically as follows. Recall that
\be 
    \Gamma_{pq}(t):=\int_S \frac{\partial g(s,t)}{\partial  s_p}N_q(s) ds, \quad 1\le p,q\le 3,
\ee
where
\be
    N=(\cos\theta \sin\phi,\sin\theta \sin\phi,\cos\phi)
\ee
and
\be
	\frac{\partial g(s,t)}{\partial  s_p} \simeq \frac{\partial g_0(s,t)}{\partial  s_p} = -\frac{s_p-t_p}{4\pi|s-t|^3}, \quad g_0(s,t):=\frac{1}{4\pi|s-t|}.
\ee
We choose a coordinate system centered at the center of the sphere such that $t=(0,0,a)$ and $s=aN$. Then
\be 
    \Gamma_{pq}(t):=-\frac{a^2}{4\pi}\int_0^{2\pi} d\theta \int_0^\pi d\phi \sin\phi \frac{(s_p-t_p)N_q}{a^3 8\sin^3\frac{\phi}{2}}, \quad 1\le p,q\le 3.
\ee
When 
\begin{align}
	&p=q=1: \quad \Gamma_{11}(t) = -1/3, \\ 
	&p=q=2: \quad \Gamma_{22}(t) = -1/3, \\ 
	&p=q=3: \quad \Gamma_{33}(t) = 1/6, \\ 
	&p \neq q: \qquad\quad \Gamma_{pq}(t) = 0.			
\end{align}
Therefore, matrix $\Gamma$ is
\be
	\Gamma\simeq\left[
	\begin{array}{ccc}
	-1/3 & 0 & 0 \\
	0 & -1/3 & 0 \\
	0 & 0 & 1/6
	\end{array}
	\right]
\ee

For example, Table \ref{tab3.1.2.0} shows the exact and asymptotic vector $Q$ when the radius of the body is $a=(1.0E-09)$ cm and the number of collocation points used to solve the integral equation \eqref{eq3.0.1} is $P=766$. Note that $a=(1.0E-09)$ cm satisfies $ka \ll 1$. The point $x_1$ in \eqref{eq1.3.8} is taken at the center of the body, the origin. Table \ref{tab3.1.2.1} and \ref{tab3.1.2.2} show the exact and asymptotic vector $E=(E_x,E_y,E_z)$, the electric field, at the point $x$ outside of the body, respectively. The distance $|x-x_1|$ is measured in cm in these tables.

\begin{table}[htbp]
  \centering
  \caption{Vector $Q_e$ and $Q_a$ when $P=766$ collocation points and $a=(1.0E-09)$ cm.}
    \begin{tabular}{rccc}
    \toprule
    \multicolumn{4}{c}{P=766, a=1.0E-09} \\
    \midrule
    \multicolumn{1}{c}{} & \multicolumn{3}{c}{ 1.0E-21 *} \\
    $Q_e$  & 0.0000 + 0.0000i & 0.0000 + 0.0000i & 0.0000 + 0.3925i \\
    $Q_a$  & 0.0000 + 0.0000i & 0.0000 + 0.0000i & 0.0000 + 0.3760i \\
    \bottomrule
    \end{tabular}%
  \label{tab3.1.2.0}%
\end{table}%
\begin{table}[htbp]
  \centering
  \caption{Vector $E_e$ for one perfectly conducting body with $a=(1.0E-09)$ cm and $P=766$ collocation points.}
    \begin{tabular}{crrr}
    \toprule
    $|x-x_1|$ & \multicolumn{3}{c}{$E_{e}(x)$} \\
    \midrule
	1.73E-08 & 1.0000 + 0.0010i & 0.0001 + 0.0000i & 0.0004 + 0.0000i \\
    1.73E-07 & 0.9999 + 0.0105i & 0.0000 + 0.0000i & 0.0000 + 0.0000i \\
    1.73E-06 & 0.9945 + 0.1045i & 0.0000 + 0.0000i & 0.0000 + 0.0000i \\	
    \bottomrule
    \end{tabular}%
  \label{tab3.1.2.1}%
\end{table}%
\begin{table}[htbp]
  \centering
  \caption{Vector $E_a$ for one perfectly conducting body with $a=(1.0E-09)$ cm and $P=766$ collocation points.}
    \begin{tabular}{crrr}
    \toprule
    $|x-x_1|$ & \multicolumn{3}{c}{$E_{a}(x)$} \\
    \midrule
	1.73E-08 & 1.0000 + 0.0010i & 0.0000 + 0.0000i & 0.0000 + 0.0000i \\
    1.73E-07 & 0.9999 + 0.0105i & 0.0000 + 0.0000i & 0.0000 + 0.0000i \\
    1.73E-06 & 0.9945 + 0.1045i & 0.0000 + 0.0000i & 0.0000 + 0.0000i \\
    \bottomrule
    \end{tabular}%
  \label{tab3.1.2.2}%
\end{table}%
\begin{table}[htbp]
  \centering
  \caption{Relative errors between the asymptotic and exact formulas for $E$ when $P=766$ collocation points and $a=(1.0E-09)$ cm.}
    \begin{tabular}{cc}
    \toprule
    $|x-x_1|$ & $E_{e}$ vs $E_{a}$ \\
    \midrule
	1.73E-08 & 4.67E-04 \\
    1.73E-07 & 4.67E-07 \\
    1.73E-06 & 4.70E-10 \\
    \bottomrule
    \end{tabular}%
  \label{tab3.1.2.3}%
\end{table}%

In this case, we also verify the following things:\\
a) Is $J$ tangential to $S$? \\
In fact, this vector $J$ is tangential to the surface $S$ of the body, $J\cdot N_s= O(10^{-14})$. \\
b) How accurate is the asymptotic formula \eqref{eq1.3.16} for $Q$?\\
We check the accuracy of the asymptotic formula for $Q$ in \eqref{eq1.3.16} by comparing it with the exact formula \eqref{eq1.3.5}, see Section \ref{sec3.1.1}, and the relative error is $4.21E-02$. The more collocation points used, the little this relative error is. \\
c) How accurate is the asymptotic formula \eqref{eq1.3.8} for $E$?\\
The accuracy of the asymptotic formula for $E$ in \eqref{eq1.3.8} can be checked by comparing it with the exact formula \eqref{eq1.3.2} at several points $x$ outside of the body, $|x-x_1|\gg a$ where $x_1$ is the center of the body, see the error analysis in Section \ref{sec3.1.1}. The relative errors are given in Table \ref{tab3.1.2.3}.

\begin{table}[htbp]
  \centering
  \caption{Relative errors of the asymptotic $E$ and $Q$ when $P=1386$ collocation points.}
    \begin{tabular}{ccccc}
    \toprule
    \multicolumn{5}{c}{$P=1386, |x-x_1|=1.73E-05$} \\
    \midrule
	$a$     			& 1.00E-07 & 1.00E-08 & 1.00E-09 & 1.00E-10 \\
    $E_{e}$ vs $E_{a}$ & 1.08E-06 & 1.08E-09 & 1.08E-12 & 1.12E-15 \\
    $Q_{e}$ vs $Q_{a}$ & 1.96E-02 & 1.96E-02 & 1.96E-02 & 1.89E-02 \\
    \bottomrule
    \end{tabular}%
  \label{tab3.1.2.4}%
\end{table}%

Table \ref{tab3.1.2.4} compares the asymptotic $Q_a$ versus exact $Q_e$ and asymptotic $E_a$ versus exact $E_e$, when $P=1386$ collocation points, $|x-x_1|=1.73E-05$ cm, and with various $a$. The errors shown in this table are relative errors, see  the error analysis in Section \ref{sec3.1.1}. As one can see from this table, the smaller the radius $a$ is, compared to the distance from the point of interest to the center of the body, the more precise the asymptotic formulas of $E$ and $Q$ are.

Furthermore, the numerical results also depend on the number of collocation points used. The more collocation points used, the more accurate the results is.

\subsection{EM wave scattering by one perfectly conducting ellipsoid body}\label{sec3.2}
In this section, we consider the EM wave scattering problem by a small perfectly conducting ellipsoid body. The method for solving the problem in this setting is the same as that of Section \ref{sec3.1} except that one needs to recompute the unit normal vector $N$.

To get the solution of this problem, one can follow the steps in Section \ref{sec3.1}. In particular, one needs to solve the linear system \eqref{eq3.0.7}-\eqref{eq3.0.9}.

To illustrate the idea, we use the same physical parameters as described in Section \ref{sec3.1.2}, except that the body now is an ellipsoid. Let $S$ be its smooth surface. The way we partition $S$ into many subdomains $S_{ij}$ is the same as the way we partition a spherical body as described in section \ref{sec3.1}. Then, the position of the collocation point in each subdomain $S_{ij}$ is defined by
\be \label{eq3.2.1}
    (x,y,z)_{ij} = (a\cos\theta_i \sin\phi_j,b\sin\theta_i \sin\phi_j,c\cos\phi_j),
\ee
where $a, b,$ and $c$ are the lengths of the semi-principal axes of the ellipsoid. The outward-pointing  normal vector $n$ to $S$ at this point is
\be \label{eq3.2.2}
    n_{ij}=n(\theta_i,\phi_j)=2\left(\frac{\cos\theta_i \sin\phi_j}{a},\frac{\sin\theta_i \sin\phi_j}{b},\frac{\cos\phi_j}{c}\right),
\ee
and the corresponding unit normal vector $N$ is
\be \label{eq3.2.3}
    N_{ij}=N(\theta_i,\phi_j)=n_{ij}/|n_{ij}|.
\ee

\begin{table}[htb]
  \centering
  \caption{Vector $E_e$ for one perfectly conducting ellipsoid body with $a=(1.0E-08)$ cm, $b=(1.0E-09)$ cm, $c=(1.0E-09)$ cm, and $P=1052$ collocation points.}
    \begin{tabular}{crrr}
    \toprule
    $|x-x_1|$ & \multicolumn{3}{c}{$E_{e}(x)$} \\
    \midrule
    	1.01E-07 & 0.9998 + 0.0010i & -0.0000 + 0.0000i & -0.0000 - 0.0000i \\
    	1.01E-06 & 0.9999 + 0.0105i & -0.0000 + 0.0000i & -0.0000 - 0.0000i \\
    	1.01E-05 & 0.9945 + 0.1045i & -0.0000 + 0.0000i & -0.0000 - 0.0000i \\
    \bottomrule
    \end{tabular}%
  \label{tab3.2.1}%
\end{table}%
\begin{table}[htb]
  \centering
  \caption{Vector $E_a$ for one perfectly conducting ellipsoid body with $a=(1.0E-08)$ cm, $b=(1.0E-09)$ cm, $c=(1.0E-09)$ cm, and $P=1052$ collocation points.}
    \begin{tabular}{crrr}
    \toprule
    $|x-x_1|$ & \multicolumn{3}{c}{$E_{a}(x)$} \\
    \midrule
	1.01E-07 & 1.0000 + 0.0010i & -0.0000 + 0.0000i & 0.0000 - 0.0000i \\
    1.01E-06 & 0.9999 + 0.0105i & -0.0000 + 0.0000i & 0.0000 - 0.0000i \\
    1.01E-05 & 0.9945 + 0.1045i & -0.0000 + 0.0000i & 0.0000 - 0.0000i \\    	
    \bottomrule
    \end{tabular}%
  \label{tab3.2.2}%
\end{table}%
For example, Tables \ref{tab3.2.1} and \ref{tab3.2.2} show the exact and asymptotic vector $E=(E_x,E_y,E_z)$, the electric field, got from solving this EM wave scattering problem with one perfectly conducting ellipsoid body, when the semi-principle axes of the body are $a=(1.0E-08)$ cm, $b=(1.0E-09)$ cm, $c=(1.0E-09)$ cm, and the number of collocation points is $P=1052$. Note that $a$, $b$, and $c$ satisfy $k\max(a,b,c) \ll 1$. The point $x_1$ in \eqref{eq1.3.8} is taken at the center of the ellipsoid body. Each row in Tables \ref{tab3.2.1} and \ref{tab3.2.2} shows the exact and asymptotic $E=(E_x,E_y,E_z)$, respectively, at the point $x$ outside of the body. The distance $|x-x_1|$ is measured in cm in these tables.

As for the case of one body, we need to verify the following things:\\
a) Is $J$ tangential to $S$? \\
In fact, this vector $J$ is tangential to the surface $S$ of the body, $J\cdot N_s= O(10^{-13})$. \\
b) Are $Q$ and $J$ correct? We check the relative error described in Section \ref{sec3.1.1}, $\text{Error} = \frac{|Q+\Gamma Q-RHS|}{|RHS|}=14\%$. The more collocation points used, the smaller this error is, for example, with $P=1762$ collocation points, this error is only 3.6\%.\\
c) How accurate is the asymptotic formula \eqref{eq1.3.8} for $E$?\\
The accuracy of the asymptotic formula for $E$ in \eqref{eq1.3.8} can be checked by comparing it with the exact formula \eqref{eq1.3.2} at several points $x$ outside of the body, $|x-x_1|\gg \max(a,b,c)$ where $x_1$ is the center of the body. The relative errors are given in Table \ref{tab3.2.3}.
\begin{table}[htb]
  \centering
  \caption{Relative errors between the asymptotic and exact formulas for $E$ when $P=1052$ collocation points, $a=(1.0E-08)$ cm, $b=(1.0E-09)$ cm, and $c=(1.0E-09)$ cm.}
    \begin{tabular}{cc}
	    \toprule
	    $|x-x_1|$ & $E_{e}$ vs $E_{a}$ \\
	    \midrule
		1.01E-07 & 1.73E-04 \\
	    1.01E-06 & 1.73E-07 \\
	    1.01E-05 & 1.73E-10 \\
	    \bottomrule
    \end{tabular}%
  \label{tab3.2.3}%
\end{table}%
\begin{table}[htb]
  \centering
  \caption{Relative errors of the asymptotic $E$ when $P=1052$ collocation points.}
        \begin{tabular}{ccccc}
	        \toprule
	        \multicolumn{5}{c}{$P=1052, |x-x_1|=1.73E-07$} \\
	        \midrule
			$a$     & 1.00E-07 & 1.00E-08 & 1.00E-09 & 1.00E-10 \\
		    $b$     & 1.00E-08 & 1.00E-09 & 1.00E-10 & 1.00E-11 \\
		    $c$     & 1.00E-08 & 1.00E-09 & 1.00E-10 & 1.00E-11 \\
		    $E_{e}$ vs $E_{a}$ & 2.65E-02 & 2.76E-05 & 2.76E-08 & 2.76E-11 \\	
	        \bottomrule
    \end{tabular}%
  \label{tab3.2.4}%
\end{table}%

Table \ref{tab3.2.4} shows the relative errors between the asymptotic $E$ versus exact $E$, when $P=1052$ collocation points, $|x-x_1|=1.73E-07$ cm, and with various semi-principle axes $a$, $b$, and $c$. As one can see from this table, the smaller the semi-principle axes are, compared to the distance from the point of interest to the center of the body, the more accurate the asymptotic formulas of $E$ is.

\subsection{EM wave scattering by one perfectly conducting cubic body}\label{sec3.3}
In this section, we consider the EM wave scattering problem by a small perfectly conducting cubic body. Again, the method for solving the problem in this setting is the same as that of Section \ref{sec3.1} except that one needs to recompute the unit normal vector $N$.

One can follow the steps outlined in Section \ref{sec3.1} to solve this problem. That means, one needs to solve the linear system \eqref{eq3.0.7}-\eqref{eq3.0.9}.

For illustration purpose, we use the same physical parameters as described in Section \ref{sec3.1.2}, except that the body now is a cube. Suppose the cube is placed in the first octant where the origin is one of its vertices, one can use the standard unit vectors in $\RRR$ as the unit normal vectors to the surfaces of the cube.

\begin{table}[htb]
  \centering
  \caption{Vector $E_e$ for one perfectly conducting body with $a=(1.0E-07)$ cm and $M=600$ collocation points.}
    \begin{tabular}{crrr}
    \toprule
    $|x-x_1|$ & \multicolumn{3}{c}{$E_e(x)$} \\
    \midrule
    		1.73E-04 & -0.5000 - 0.8660i & -0.0000 + 0.0000i & -0.0000 + 0.0000i \\
    	    1.73E-05 & 0.5000 + 0.8660i & 0.0000 + 0.0000i & 0.0000 + 0.0000i \\
    	    1.73E-06 & 0.9945 + 0.1045i & 0.0006 + 0.0000i & 0.0006 + 0.0000i \\
    \bottomrule
    \end{tabular}%
  \label{tab3.3.1}%
\end{table}%
\begin{table}[htb]
  \centering
  \caption{Vector $E_a$ for one perfectly conducting body with $a=(1.0E-07)$ cm and $M=600$ collocation points.}
    \begin{tabular}{crrr}
    \toprule
    $|x-x_1|$ & \multicolumn{3}{c}{$E_a(x)$} \\
    \midrule
    		1.73E-04 & -0.5000 - 0.8660i & 0.0000 - 0.0000i & 0.0000 + 0.0000i \\
    	    1.73E-05 & 0.5000 + 0.8660i & -0.0000 + 0.0000i  & -0.0000 - 0.0000i \\
    	    1.73E-06 & 0.9945 + 0.1045i & -0.0000 + 0.0000i  & -0.0000 - 0.0000i \\    	
    \bottomrule
    \end{tabular}%
  \label{tab3.3.2}%
\end{table}%
For example, Tables \ref{tab3.3.1} and \ref{tab3.3.2} show the exact and asymptotic vector $E=(E_x,E_y,E_z)$, the electric field, got from solving this EM wave scattering problem with one perfectly conducting cubic body, when the half side of the body is $a=(1.0E-07)$ cm and the number of collocation points is $M=600$. Note that $a=1.0E-07$ cm satisfies $ka \ll 1$. The point $x_1$ in \eqref{eq1.3.8} is taken at the center of the cubic body. Each row in Tables \ref{tab3.3.1} and \ref{tab3.3.2} shows the exact and asymptotic $E=(E_x,E_y,E_z)$, respectively, at the point $x$ outside of the body. The distance $|x-x_1|$ is measured in cm in these tables.

As before, for the case of one body, we need to verify the following things:\\
a) Is $J$ tangential to $S$? \\
In fact, this vector $J$ is tangential to the surface $S$ of the body, $J\cdot N_s= O(10^{-13})$. \\
b) How accurate is the asymptotic formula \eqref{eq1.3.16} for $Q$?\\
We check the relative error described in Section \ref{sec3.1.1}, $\text{Error} = \frac{|Q+\Gamma Q-RHS|}{|RHS|}=1.13\%$. \\
c) How accurate is the asymptotic formula \eqref{eq1.3.8} for $E$?\\
The accuracy of the asymptotic formula for $E$ in \eqref{eq1.3.8} can be checked by comparing it with the exact formula \eqref{eq1.3.2} at several points $x$ outside of the body, $|x-x_1|\gg a$ where $x_1$ is the center of the body. The relative errors are given in Table \ref{tab3.3.3}.
\begin{table}[htb]
  \centering
  \caption{Relative errors between the asymptotic and exact formulas for $E$ when $M=600$ collocation points and $a=(1.0E-07)$ cm.}
    \begin{tabular}{cc}
	    \toprule
	    $|x-x_1|$ & $E_e$ vs $E_a$ \\
	    \midrule
			1.73E-03 & 1.19E-08 \\
		    1.73E-04 & 1.19E-07 \\
		    1.73E-05 & 1.52E-06 \\
		    1.73E-06 & 8.64E-04 \\		
	    \bottomrule
    \end{tabular}%
  \label{tab3.3.3}%
\end{table}%

\begin{table}[htb]
  \centering
  \caption{Relative errors of the asymptotic $E$ when $M=600$ collocation points.}
        \begin{tabular}{cccc}
	        \toprule
	        \multicolumn{4}{c}{$M=600, |x-x_1|=1.73E-06$} \\
	        \midrule
	      		 $a$     & 1.00E-07 & 1.00E-08 & 1.00E-09 \\
	             $E_e$ vs $E_a$ & 8.64E-04 & 6.49E-07 & 6.32E-10 \\	
	        \bottomrule
    \end{tabular}%
  \label{tab3.3.4}%
\end{table}%

Table \ref{tab3.3.4} shows the relative errors between the asymptotic $E$ versus exact $E$, when $M=600$ collocation points, $|x-x_1|=(1.73E-06)$ cm, and with various $a$. From this table, we can see that the smaller the side of the cube is, compared to the distance from the point of interest to the center of the body, the more accurate the asymptotic formula of $E$ is.

\subsection{EM wave scattering by many small perfectly conducting bodies}
To illustrate the idea, consider a domain $\Omega$ as a unit cube placed in the first octant such that the origin is one of its vertices. This domain $\Omega$ contains $M$ small bodies. Suppose these small bodies are particles. We use GMRES iterative method, see \cite{GMRES}, to solve the linear system \eqref{eq3.5.7}-\eqref{eq3.5.9}. The following physical parameters are used to solve the EM wave scattering problem
\begin{itemize}
     \item Speed of wave, $c=(3.0E+10)$ cm/sec.
     \item Frequency, $\omega=(5.0E+14)$ Hz.
     \item Wave number, $k =(1.05E+05)$ cm$^{-1}$.
     \item Wave length, $\lambda= (6.00E-05)$ cm.
     \item Direction of incident plane wave, $\alpha = (0,  1,  0)$.
     \item Magnetic permeability, $\mu = 1$.
     \item Volume of the domain $\Omega$ that contains all the particles, $|\Omega| = 1$ cm$^3$.
     \item The distance between two neighboring particles, $d = (1.00E-07)$ cm.
     \item Vector $\mathcal{E} = (1,  0,  0)$.
     \item Vector $A_0$: $A_{0m}:=(I+\Gamma)^{-1}\nabla \times E_0(x)|_{x=x_m}=(I+\Gamma)^{-1}\nabla \times \mathcal{E} e^{ik\alpha\cdot x}|_{x=x_m}$.
\end{itemize}
Note that the distance $d$ satisfies the assumption $d \ll \lambda$. The radius $a$ of the particles is chosen variously so that it satisfies the assumption $ka \ll 1$. For illustration purpose, the problem of EM wave scattering by many small perfectly conducting bodies is solved with $M=27$ and $1000$ particles.

\begin{table}[htbp]
  \centering
  \caption{Vector $E$ when $M=27$ particles, $d=(1.0E-07)$ cm and $a=(1.0E-09)$ cm.}
    \begin{tabular}{rrr}
       \toprule
       		\multicolumn{3}{c}{$M=27$, $d=1.0E-07$, $a=1.0E-09$}\\
       \midrule
			1.00E+00+1.01E-14i &  5.69E-17-1.01E-14i &  0.00E+00+0.00E+00i \\
		    1.00E+00+1.19E-14i &  0.00E+00+0.00E+00i &  0.00E+00+0.00E+00i \\
		    1.00E+00+1.01E-14i &  -5.69E-17+1.01E-14i &  0.00E+00+0.00E+00i \\
		    1.00E+00+1.05E-02i &  1.24E-16-1.19E-14i &  -1.36E-29-5.20E-36i \\
		    1.00E+00+1.05E-02i &  0.00E+00+0.00E+00i &  -4.80E-30-5.20E-36i \\
		    1.00E+00+1.05E-02i &  -1.24E-16+1.19E-14i &  -1.22E-30-5.20E-36i \\
		    1.00E+00+2.09E-02i &  1.54E-16-1.01E-14i &  -3.40E-30-1.04E-35i \\
		    1.00E+00+2.09E-02i &  0.00E+00+0.00E+00i &  -2.43E-30-1.04E-35i \\
		    1.00E+00+2.09E-02i &  -1.54E-16+1.01E-14i &  -1.20E-30-1.04E-35i \\
		    1.00E+00+1.19E-14i &  6.63E-17-1.19E-14i &  0.00E+00+0.00E+00i \\
		    1.00E+00+1.40E-14i &  4.80E-30+0.00E+00i &  0.00E+00+0.00E+00i \\
		    1.00E+00+1.19E-14i &  -6.63E-17+1.19E-14i &  0.00E+00+0.00E+00i \\
		    1.00E+00+1.05E-02i &  1.47E-16-1.40E-14i &  -4.80E-30-5.20E-36i \\
		    1.00E+00+1.05E-02i &  2.61E-30+0.00E+00i &  -2.61E-30-5.20E-36i \\
		    1.00E+00+1.05E-02i &  -1.47E-16+1.40E-14i &  -9.24E-31-5.20E-36i \\
		    1.00E+00+2.09E-02i &  1.82E-16-1.19E-14i &  -2.43E-30-1.04E-35i \\
		    1.00E+00+2.09E-02i &  9.24E-31+0.00E+00i &  -1.85E-30-1.04E-35i \\
		    1.00E+00+2.09E-02i &  -1.82E-16+1.19E-14i &  -1.01E-30-1.04E-35i \\
		    1.00E+00+1.01E-14i &  5.69E-17-1.01E-14i &  0.00E+00+0.00E+00i \\
		    1.00E+00+1.19E-14i &  2.43E-30+0.00E+00i &  0.00E+00+0.00E+00i \\
		    1.00E+00+1.01E-14i &  -5.69E-17+1.01E-14i &  0.00E+00+0.00E+00i \\
		    1.00E+00+1.05E-02i &  1.24E-16-1.19E-14i &  -1.22E-30-5.20E-36i \\
		    1.00E+00+1.05E-02i &  1.85E-30+0.00E+00i &  -9.24E-31-5.20E-36i \\
		    1.00E+00+1.05E-02i &  -1.24E-16+1.19E-14i &  -5.03E-31-5.20E-36i \\
		    1.00E+00+2.09E-02i &  1.54E-16-1.01E-14i &  -1.20E-30-1.04E-35i \\
		    1.00E+00+2.09E-02i &  9.98E-31+0.00E+00i &  -1.01E-30-1.04E-35i \\
		    1.00E+00+2.09E-02i &  -1.54E-16+1.01E-14i &  -6.54E-31-1.04E-35i \\
        \bottomrule
    \end{tabular}%
  \label{tab1.0}%
\end{table}%

For example, Tables \ref{tab1.0} show the result of solving the EM wave scattering problem with $M=27$ particles in the unit cube in which the distance between neighboring particles is $d=(1.0E-07)$ cm and the radius of the particles is $a=(1.0E-09)$ cm. Each row in Tables \ref{tab1.0} is a vector $E(i)=(E_x,E_y,E_z)(i)$ at the point $i$ in the cube. The norm of this asymptotic solution $E$ is $5.20E+00$ and the error of the solution is $8.16E-10$. This error is computed using \eqref{eq3.5.51}.

\begin{figure}[htbp]
    \centering
    \includegraphics[scale=0.95]{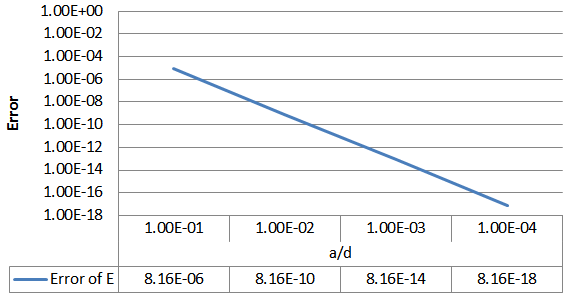}
    \caption{Error of the asymptotic solution $E$ when $M=27$ and $d=(1.0E-07)$ cm.}
    \label{fig1}
\end{figure}

Table \ref{tab1} shows the relative errors of $E$ when there are $M=27$ particles in the cube, the distance between neighboring particles is $d=(1.0E-07)$ cm, and with various radius $a$. Figure \ref{fig1} shows the relative error of the asymptotic $E$. From this figure, one can see that when the ratio $a/d$ decreases from $1.0E-01$ to $1.0E-04$, the error of the asymptotic solution decreases linearly and rapidly from $8.16E-06$ to about $8.16E-18$. The smaller the ratio $a/d$ is, the better the asymptotic formula \eqref{eq2.1.22} approximates $E$.

\begin{table}[htbp]
  \centering
  \caption{Error of the asymptotic solution $E$ when $M=27$ and $d=(1.0E-07)$ cm.}
    \begin{tabular}{ccccc}
      \toprule
        \multicolumn{5}{c}{M=27, d=1.0E-07} \\
      \midrule
		a     & 1.00E-08 & 1.00E-09 & 1.00E-10 & 1.00E-11 \\
	    a/d   & 1.00E-01 & 1.00E-02 & 1.00E-03 & 1.00E-04 \\
	    Norm of E & 5.20E+00 & 5.20E+00 & 5.20E+00 & 5.20E+00 \\
	    Error of E & 8.16E-06 & 8.16E-10 & 8.16E-14 & 8.16E-18 \\
      \bottomrule
    \end{tabular}%
  \label{tab1}%
\end{table}%

\begin{figure}[htbp]
    \centering
    \includegraphics[scale=0.95]{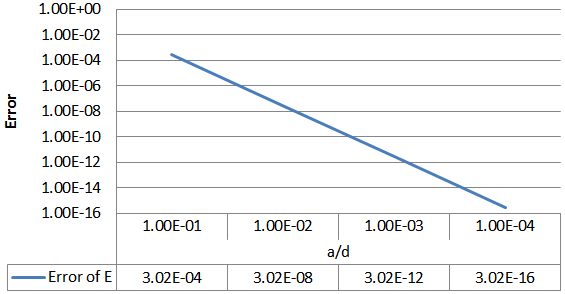}
    \caption{Error of the asymptotic solution $E$ when $M=1000$ and $d=(1.0E-07)$ cm.}
    \label{fig2}
\end{figure}

Table \ref{tab2} and Figure \ref{fig2} show the results of solving the problem with $M=1000$ particles, when the distance between neighboring particles is $d=(1.0E-07)$ cm, and with different radius $a$. From these table and figure, one can see that the relative error of the asymptotic solution in this case is also very small, less than $3.02E-04$, when the ratio $a/d < 1.0E-01$. In this case, the error of the asymptotic $E$ is greater than that of the previous case when $M=27$. However, this time, the error is also decreasing quickly and linearly when the ratio $a/d$ decreases from $1.0E-01$ to $1.0E-04$. Therefore, the asymptotic formula \eqref{eq2.1.22} for the solution $E$ is applicable when $a \ll d$.

\begin{table}[htbp]
  \centering
  \caption{Error of the asymptotic solution $E$ when $M=1000$ and $d=(1.0E-07)$ cm.}
    \begin{tabular}{ccccc}
      \toprule
         \multicolumn{5}{c}{M=1000, d=1.0E-07} \\
      \midrule
		a     & 1.00E-08 & 1.00E-09 & 1.00E-10 & 1.00E-11 \\	
	    a/d   & 1.00E-01 & 1.00E-02 & 1.00E-03 & 1.00E-04 \\
	    Norm of E & 3.16E+01 & 3.16E+01 & 3.16E+01 & 3.16E+01 \\
	    Error of E & 3.02E-04 & 3.02E-08 & 3.02E-12 & 3.02E-16 \\
      \bottomrule
    \end{tabular}%
  \label{tab2}%
\end{table}%

\section{Conclusions} \label{sec4}
In this paper, we present a numerical method for solving the EM wave scattering by one and many small perfectly conducting bodies. One of the advantages of this method is that it is relatively easy to implement. Furthermore, one can get an asymptotically exact solution to the problem when the characteristic size of the bodies tends to zero. To illustrate the applicability and efficiency of the method, we use it to solve the EM wave scattering problem by one and many small perfectly conducting bodies. Numerical results of these experiments are presented and error analysis of the asymptotic solutions for the case of one and many bodies are also discussed. For the case of one small body, one can always find the exact solution using the described method. For the case of many small bodies, the accuracy of our method is high if $a \ll d \ll \lambda$. 

The problem of EM wave scattering is much harder to treat, compared to scalar wave scattering \cite{Nakayama1981, Ito1985, TMaterial, TFastScalar}. For scalar wave scattering problem, a fast algorithm is developed in \cite{TFastScalar} to deal with billions of particles. This is still open to EM wave scattering. One might consider this as future research.

\bibliographystyle{ieeetran}

\end{document}